\documentclass[final]{siamltex}

\usepackage{psfig}       
\usepackage{graphicx}   
\usepackage{amsfonts}  
\newtheorem{thm}{Theorem}[section] 

\newtheorem{cor}[thm]{Corollary}
\newtheorem{defn}[thm]{Definition}

\renewcommand{\theequation}{\thesection.\arabic{equation}}
\newcommand{\field}[1]{\mathbb{#1}} 

\def\z{{\bar{z}}}

\def\P{\tilde{P}}

\def\R{\mathbb{R}}

\def\th{\theta}

\title{On Numerical Algorithms for the Solution of a Beltrami Equation}

\author{Denis Gaydashev\thanks{Department of Mathematics, University of Toronto,
Toronto, Ontario, Canada M5S 3G3 ({\tt gaidash@math.toronto.edu}).}
	\and \framebox{Dmitry Khmelev} \thanks{This project was completed by the first author after D. Khmelev's sudden passage.}}

\begin{document}

\maketitle

\begin{abstract}
The paper concerns numerical algorithms for solving the Beltrami equation $f_{\bar{z}}(z)=\mu(z)  f_z(z)$ for a compactly supported $\mu$.

First, we study an efficient algorithm that has been proposed in the literature, and  present its rigorous justification.  We then propose a different scheme for solving the Beltrami equation which has  a comparable speed and accuracy, but has the virtue of a greater simplicity of implementation. 
\end{abstract}



\pagestyle{myheadings}
\thispagestyle{plain}
\markboth{D. GAYDASHEV AND D. KHMELEV}{ON NUMERICAL ALGORITHMS FOR BELTRAMI EQUATION}


\setcounter{page}{1}
\section{Introduction}\label{intro}

The of existence of a quasiconformal mapping $g$ with a given complex dilatation $\mu$ is one of the central problems of classical complex analysis.

We will first recall a  definition of a quasiconformal map. Given an open set  $\Omega \subset \field{C}$, a map $f: \Omega \rightarrow f(\Omega) \subset \field{C}$ is said to be absolutely continuous on lines (ACL) in a rectangle $R \subset \Omega$ with sides parallel to the $x$ and $y$ axes if $f$ is absolutely continuous on almost every horizontal and vertical line in $R$. The map $f$ is ACL on $\Omega$, if it is ACL on every rectangle in $\Omega$. Partial derivatives $f_z=(f_x-i f_y)/2$ and $f_{\bar{z}}=(f_x+i f_y)/2$ of such map exist a.e.  in $\Omega$.
\begin{defn}
A homeomorphism  $f: \Omega \rightarrow f(\Omega)$  is K-quasiconformal if and only if the following holds:
\begin{itemize}
\item[i)] $f$ is ACL on every rectangle in $\Omega$,
\item[ii)] $|f_{\bar{z}}| \le {K-1 \over K+1} |f_z|$ a.e. in $\Omega$.
\end{itemize}
\end{defn}

The complex dilatation of a quasiconformal homeomorphism $f$ defined a.e. is 
\begin{eqnarray}
\nonumber \mu_f(z)&=&f_{\bar{z}}(z)/f_z(z).
\end{eqnarray}

  A mapping $f$ with a prescribed complex dilatation $\mu$ solves the following {\it Beltrami equation}: 
\begin{equation}\label{belt_eqn}
f_{\bar{z}}(z)=\mu(z)  f_z(z).
\end{equation}

The question of existence of a solution of  $(\ref{belt_eqn})$ is addressed in the famous Ahlfors--Bers--Boyarskii theorem (see \cite{AB} and \cite{Boyarskii}):
\begin{thm} \label{ABB}(Ahlfors and Bers, Boyarskii).
Let $\mu \in L_\infty(\hat{\field{C}})$ satisfy  $\left\| \mu \right\|_\infty \le K < 1$. Then there exists a unique map $f^\mu: \hat{\field{C}} \rightarrow \hat{\field{C}}$ such that $f^\mu$ satisfies the Beltrami equation $f^\mu_{\bar{z}} (z)=\mu(z) f^\mu_z(z)$, $f^\mu$ fixes $0$, $1$ and $\infty$, and $f^\mu$ is a $(1+K)/(1-K)$-quasiconformal map. 
\end{thm}

According to the classical proof of this Theorem, a solution of the Beltrami equation can be constructed with the help of two operators: the first operator being the Hilbert transform
\begin{equation}\label{Hilbert_tran}
T [h](z)={i \over 2 \pi} \lim_{\epsilon \to 0} \int \int_{\field{C} \setminus B(z,\epsilon)} {h(\xi) \over (\xi -z)^2} \ d \bar\xi \wedge d \xi, 
\end{equation}
the second --- the Cauchy transform
\begin{equation}\label{Cauchy_tran}
P [h] (z)={i \over 2 \pi} \int \int_{\field{C}} {h(\xi) \over \xi -z} - {h(\xi) \over \xi} \ d \bar{\xi} \wedge d \xi.
\end{equation}

If the differential $\mu$ is compactly supported, then the quasiconformal homeomorphism $f$ is constructed in two steps. First, one solves the following equation:
\begin{equation}  \label{eq:h=Tmuh+Th}
  h^*=T[\mu h^*]+T[\mu].
\end{equation}

If $h^*$ satisfying (\ref{eq:h=Tmuh+Th}) is found, then
$f$ solving (\ref{belt_eqn}) is given by
\begin{eqnarray}
\nonumber   f(z)&=&P[\mu(h^*+1)](z)+z.
\end{eqnarray}

The operator $P$ is well-defined and is H\"older continuous for $h \in L^{p}(\field{C})$, $p>2$. This implies that a solution to (\ref{eq:h=Tmuh+Th}) must also
be constructed in $L_p(\field{C})$, $p>2$. The celebrated Calderon--Zygmund inequality (cf \cite{CalZyg}) states that
\begin{eqnarray}
\nonumber   \|T[h]\|_p&\le&C_p\|h\|_p
\end{eqnarray}
for all $h\in L_p(\field{C})$ where $C_p$ is some constant, depending on $p$
only, which can be chosen so that $C_p\to1$ as $p\to 2$. While the
optimal value of $C_p$ for all $p$ is not known, it is known that
$C_2=1$, and $T$ is an isometry in $L_2(\field{C})$: 
\begin{eqnarray}
\nonumber \|T[h]\|_2&=&\|h\|_2.
\end{eqnarray}

To solve (\ref{eq:h=Tmuh+Th}) one has to choose a $p$ such that
$K C_p<1$. Then the operator $h\mapsto{}T[\mu{}h]$ is  a contraction in $L_p(\field{C})$ and iterations
\begin{equation}\label{Hilbert_iter}
h^{n+1}=T[\mu{}h^n]+T[\mu]
\end{equation}
converge to a solution of (\ref{eq:h=Tmuh+Th}): $h^n\to{}h^*$
as $n\to\infty$.

We can now state the extended version of the Ahlfors--Bers--Boyarskii theorem  for compactly supported differentials:
\begin{thm} \label{compactABB}
Let $\mu \in L_\infty(\field{C})$ be compactly supported and satisfy $\left\| \mu \right\|_\infty \le K < 1$. 
Then for every $p>2$ such that $C_p K<1$, the operator $h \rightarrow T [ \mu (h+1)] $ is a contraction on $L_p(\field{C})$ with a unique fixed point $h^*$. Moreover, the solution of the Beltrami equation, $f^\mu$, is  given by 
\begin{equation}\label{solution}
 f^\mu=P [\mu (h^*+1)] +id,
\end{equation}
and is such that $f(0)=0$, $f$ is continuous, has distributional derivatives, and $f_z-1 \in L_p(\field{C})$.
\end{thm}

We will concentrate on the numerical implementation of the procedure described in this Theorem. This is an important practical problem which has been addressed by various authors. 

Efficient algorithms for computing the two required transforms $T$ and $P$ have been proposed in \cite{Daripa3}, \cite{Daripa1} and \cite{Daripa2}. We will proceed with a very brief description of these algorithms.

As we have pointed out, the Cauchy transform of an $L_p$-function, $p > 2$, is well-defined and is known to be H\"older continuous with exponent $1-2/p$ (cf \cite{Ahlfors}, \cite{Markovic}). In \cite{Daripa3} an effective algorithm for the numerical computation of the following integral transform was suggested:
\begin{equation}  \label{eq:Pt}
  \P [h] (z)={i \over 2 \pi} \int \int_{\field{C}} {h(\xi) \over (\xi -z)} \ d \bar{\xi} \wedge d \xi.
\end{equation}

A comparison with (\ref{Cauchy_tran}) immediately shows that,
\begin{eqnarray} \label{P_Pt}
  P [h](z)&=&\P [ h ](z)-\P [ h ] (0).
\end{eqnarray}

 Represent $h$ and $\P [h]$ as a  Fourier series:
\begin{eqnarray}
\label{h_coeff} h(r e^{i \theta})&=&\sum^{\infty}_{k=-\infty} h_k(r) e^{i k \theta}, \\
\label{Pt_coeff} \P [ h ](r e^{i \theta})&=&\sum^{\infty}_{k=-\infty} p_k(r) e^{i k \theta},
\end{eqnarray}
where the coefficients of the $\P$-transform are given by 
\begin{equation} \label{pk}
p_k(r)= {1 \over 2 \pi } \int^{2 \pi}_0 e^{-i k \theta} \P[h](r e^{i \theta}) \ d \theta,
\end{equation}
then, according to \cite{Daripa3}, 
\begin{equation}\label{Cauchy_alg}
 p_k(r)=\left\{
\begin{array}{cc}
  \displaystyle    2\int_0^r\left(\frac{r}{\rho}\right)^k h_{k+1}(\rho)d\rho,&k < 0,\\
  \displaystyle   -2\int_r^\infty\left(\frac{r}{\rho}\right)^k h_{k+1}(\rho)d\rho,&k\ge0.
\end{array}
 \right.
\end{equation}

Next, to obtain similar formulae for the Hilbert transform, assume that $h$ is a H\"older continuous function compactly supported in an open  disk around zero of radius $R$, $B(0,R) \subset \field{C}$.  The Hilbert transform of such function is well-defined (see next Section). Represent the Hilbert transform as a Fourier series
\begin{eqnarray}\label{T_coeff} 
T [ h ](r e^{i \theta})&=&\sum^{\infty}_{k=-\infty} c_k(r) e^{i k \theta}, \quad c_k(r)= {1 \over 2 \pi } \int^{2 \pi}_0 e^{-i k \theta} T[h](r e^{i \theta}) \ d \theta. 
\end{eqnarray}

In \cite{Daripa1} and \cite{Daripa2} the authors arrive at the following expressions for these coefficients:
\begin{eqnarray}
\label{Hilbert_alg_1}c_0(0)&=& -2 \lim_{\epsilon \rightarrow 0} \int^R_\epsilon {h_2(\rho) \over \rho}  d \rho, \ {\rm and} \ c_k(0)=0, \ {\rm whenever} \ k \ne 0, \\
\label{Hilbert_alg_2}c_k(r)&=&A_k \int^r_0 { r^k  \over \rho^{k+1} }  h_{k+2}(\rho)   d \rho+B_k\int^R_r  { r^k  \over \rho^{k+1} }  h_{k+2}(\rho)   d \rho +h_{k+2}(r),
\end{eqnarray}
where 
\begin{eqnarray} \label{A_k_B_k} 
A_k&=&\left\{ 0 \ , \ k \ge 0,  \atop  2 (k+1) \ , k<0, \right.   \ {\rm and} \    B_k=\left\{-2  (k+1),  \ k \ge 0, \atop 0, \ k<0. \right.
\end{eqnarray}

Given  values of $h$, for instance, on a circular $N \times M$ grid that contains the compact support of $h$, one can use a fast Fourier transform (FFT, cf \cite{NR}) to find the values of the coefficients $h_k$ at the radii $r_i$, $ 1 \le i\le M$. Next, one can use these values to construct a piecewise constant, a piecewise linear or a spline approximation of the functions $h_k$ (the choice of approximation, of course, depends on the known or expected smoothness of $h_k$). This allows one to compute integrals in (\ref{Cauchy_alg}) and in (\ref{Hilbert_alg_1})--(\ref{Hilbert_alg_2}). Armed with these implementations of the Hilbert and Cauchy transforms, one can try to solve the Beltrami equation (\ref{belt_eqn}), first by running iterations (\ref{Hilbert_iter}) for some time, and, finally, applying (\ref{solution}). It is convenient to use the point-wise multiplication of grid values of $h$ and $\mu+1$ inside the Hilbert transform in (\ref{Hilbert_iter}), rather than the multiplication of their Fourier series: The order of the computational complexity of the point-wise multiplication is $O(NM)$, as opposed to $O(NM^2)$ for the series. The transition from the representation of $h$ as a Fourier series to point values at each iteration step can be performed with the help of the  FFT. This way, the computational complexity of one iteration step becomes $O(N M \log_2 M)$.

 The objective of our paper will be two-fold. First, we will present a rigorous mathematical justification of the formulae (\ref{Hilbert_alg_1})--(\ref{Hilbert_alg_2}) for the singular operator $T$ different from that of \cite{Daripa2}. Second, we will propose an algorithm which avoids programming formulae (\ref{Hilbert_alg_1})--(\ref{Hilbert_alg_2}) altogether. It should be noted, that in our experience writing the code for the Hilbert transform can be rather cumbersome, therefore this second algorithm can be viewed as an improvement.

\medskip\section{A Justification of P. Daripa and D. Mashat's algorithm for the Hilbert transform}\label{sect1}
\setcounter{equation}{0}

The Hilbert transform is a singular integral transform, therefore, first of all, one has to specify the conditions under which this integral can be computed. The issue of existence of a rather general class of singular integral transforms (including the Hilbert transform) has been addressed in the classical work \cite{CalZyg} of A. P. Calderon and A. Zygmund. The sufficient condition for the existence of the Hilbert transform itself can be stated in a very concise way (cf. \cite{Ahlfors}, \cite{Carleson}):

\begin{lemma}\label{exist_lemma}
If $h$ is H\"older continuous, then the Hilbert integral transform of $h$ exists as a Cauchy principle value. 
\end{lemma}

In this light, the main result of \cite{Daripa2} can be stated as follows:

\begin{thm}\label{Daripa_thm}
If $T[h](z)$ exists  as a Cauchy principle value, than the Fourier coefficients of this Hilbert transform are given by formulae (\ref{Hilbert_alg_1})--(\ref{Hilbert_alg_2}).
\end{thm}

The proof of this Theorem presented in \cite{Daripa2} starts with an observation that according to (\ref{Hilbert_tran}) and (\ref{T_coeff}), coefficients of the Hilbert transform are given by the following formula:
\begin{eqnarray}\label{ck_long}
 c_k(r)=  -{1 \over 2 \pi^2 } \int^{2 \pi}_0 e^{-i k \theta}  \lim_{\epsilon \to 0} \int \int_{\field{C} \setminus B(r e^{i \theta},\epsilon)} {h(\rho e^{i \alpha}) \over (\rho e^{i \alpha} -r e^{i \theta} )^2} \ \rho d \rho d \alpha d \theta. 
\end{eqnarray}

P. Daripa and D. Mashat proceed in their proof by computing the Hilbert transform in the problematic region ($\rho$ close to $r$) explicitly.  A crucial step in the ensuing proof is the interchange of the order of integration in (\ref{ck_long}). We feel that the possibility of this interchange has to be carefully explained: Although the argument that we present in this Section is not exceptionally difficult, it nevertheless involves some reasoning and has to be carefully carried out.

We will now proceed with our argument. Let $h$ be H\"older continuous  with exponent $\gamma$ and constant $A$, and let it be  compactly supported in an open  disk around zero of radius $R$, $B(0,R) \subset \field{C}$. Notice, that the variable $z$ in (\ref{Hilbert_tran}) can range over all of $\field{C}$ and not only over $B(0,R)$. We will start by showing that the Hilbert transform of a  H\"older continous function is bounded.

Choose a small $\delta<1/R$. For all $\epsilon < \delta R$  define 
\begin{eqnarray}\label{eq1}
\nonumber T_\epsilon[h] (z)&=&{i \over 2 \pi} \int   \int_{B(0,R) \setminus B(z,\epsilon)}{ h(\xi) \over (\xi -z)^2 }\ d \bar\xi \wedge d \xi \\
\nonumber &=&{i \over 2 \pi} \int \int_{B(0,R) \setminus B(z,\delta R)} {h(\xi) \over (\xi -z)^2} \ d \bar\xi \wedge d \xi \\ 
\nonumber &+& {i \over 2 \pi} \int \int_{B(z,\delta R) \setminus B(z,\epsilon)} { h(\xi) \over (\xi -z)^2 } \ d \bar{\xi} \wedge d \xi\\ 
&=&T^1[h](z)+T^2_\epsilon [h](z), 
 \end{eqnarray} 
 \noindent where $r=|z|$. 
\begin{lemma}\label{lemma1}
There exists a constant $C$ such that $\left| T^2_\epsilon [h](z) \right| < C$ for all $0<\epsilon < \delta R$ and all $z \in \overline{ B(0,R) }$. 
\end{lemma}
\begin{proof}
\begin{eqnarray}
\nonumber \left| T^2_\epsilon [h](z) \right|&=&\left| {i \over 2 \pi} \int  \int_{B(z,\delta R) \setminus B(z,\epsilon)}  { h(\xi)  \over  (\xi -z)^2} \  d \bar{\xi} \wedge d \xi \right| \\
\nonumber &=&\left|-{ 1 \over \pi}  \int^{ \delta R}_\epsilon \int^{2 \pi}_{0}   {h(z+\rho e^{i \theta}) \over \rho^2 e^{2 i \theta }}  \rho  \ d \rho d \theta \right|\\
\nonumber &=&\left|-{ 1 \over \pi}   \int^{\delta R}_\epsilon \int^{\pi \over 2}_0  { h  \left(z+\rho e^{i \theta } \right) - h \left( z+\rho e^{i (\theta +\pi / 2 ) } \right)  \over \rho  e^{2 i \theta } }   d  \rho d \theta \right.\\
\nonumber &-& \left. { 1 \over \pi}   \int^{\delta R}_\epsilon  \int^{\pi \over 2}_0  {h  \left( z + \rho e^{i (\theta + \pi) } \right) - h  \left(z + \rho e^{i ( \theta + 3 \pi /  2)} \right)   \over  \rho  e^{2 i \theta } }   d \rho d \theta \right| \\
\nonumber &\le& { 1 \ \over \pi}  \int^{\delta R}_\epsilon   \int^{\pi \over 2}_0 {\left| h  \left(z + \rho e^{i \theta} \right) - h \left(z + \rho e^{i (\theta + \pi /  2) }\right) \right| \over \rho} \ d \rho d \theta  \\
\nonumber &+& { 1 \over \pi}  \int^{\delta R}_\epsilon \int^{\pi \over 2}_0 {\left| h \left(z + \rho e^{i (\theta + \pi)} \right) - h  \left(z + \rho e^{i (\theta + 3 \pi / 2 )} \right) \right| \over \rho} \ d \rho d \theta \\
\nonumber & \le & {1 \over 2}  \int^{\delta R}_\epsilon {A \rho^\gamma \left|1-e^{i {\pi \over 2}} \right|^\gamma \over \rho}  \ d \rho + {1 \over 2}  \int^{\delta R}_\epsilon {A \rho^\gamma \left|1-e^{i {\pi \over 2}}\right|^\gamma \over \rho} \ d \rho \\
\nonumber & \le &  A  \left|1-i \right|^\gamma {1 \over \gamma} \left( \delta R \right)^\gamma=C, 
\end{eqnarray}
\noindent where we have used H\"older continuity.
\end{proof} \bigskip

Next, we will derive the expressions for the Fourier coefficients of the Hilbert transform separately for  $z \in {\rm int}  B(0,R)  -\{0\}$,  $z \in \partial B(0,R)$, $z=0$ and $z \in \field{C} \setminus  \overline{B (0,R) }$.

\subsection{Case $z \in {\rm int} B(0,R)-\{0\}$}\label{sect2}
$\phantom{aaaa}$
\bigskip

\begin{cor} \label{Lebesgue}
\begin{eqnarray}
\nonumber c_k(r)&=&-{1 \over 2 \pi^2}  \lim_{\epsilon \to 0} \int^{r-\epsilon}_0  \int^{2 \pi}_0  h(\rho e^{i \alpha}) \int^{2 \pi}_0  { e^{-i k \theta} \over (\rho e^{ i \alpha} -r e^{i \theta})^2} \  d \theta d \alpha  \rho d \rho\\
\nonumber &-&{1 \over 2 \pi^2} \lim_{\epsilon \to 0} \int^{R}_{r+\epsilon}  \int^{2 \pi}_0  h(\rho e^{i \alpha}) \int^{2 \pi}_0  { e^{-i k \theta} \over (\rho e^{ i \alpha} -r e^{i \theta})^2} \  d \theta d \alpha  \rho d \rho +h_{k+2}(r)
\end{eqnarray}
\noindent for all $0<r <R$.
\end{cor}
\begin{proof} 
 Let $z=r e^{i \theta}$. Recall, that $T_\epsilon[h]=T^1[h]+T^2_\epsilon [h]$ where $T^1[h]$ and $T^2_\epsilon [h]$ are as in (\ref{eq1}). The first of these operator is clearly bounded, the second is bounded by the preceding Lemma. Therefore,  $|T_\epsilon[h](z)|$ is bounded by a constant function, while the function $\psi_\epsilon(\theta)= e^{-i \theta} T_\epsilon[h]( r  e^{i \theta} )$ is bounded by ${\rm const} \cdot  e^{-i \theta}$ which is certainly in $L_1([0,2 \pi))$. A direct application of the Lebesgue Dominated Convergence Theorem shows that the limiting procedure and the integration over the circle can be interchanged:
\begin{eqnarray}
\nonumber c_k(r)&=&{i \over 4 \pi^2}\int^{2 \pi}_0 \lim_{\epsilon \to 0}  \psi_\epsilon (\theta) d \theta={i \over 4 \pi^2}\lim_{\epsilon \to 0} \int^{2 \pi}_0  \psi_\epsilon (\theta) d \theta \\
\nonumber & =&{i \over 4 \pi^2} \lim_{\epsilon \to 0} \int^{2 \pi}_0 \int  \int_{B(0,R) \setminus B(z,\epsilon)} h(\xi)  { e^{-i k \theta} \over (\xi -r e^{i \theta})^2} \  d \bar{\xi} \wedge d \xi d \theta= \lim_{\epsilon \to 0}  c^\epsilon_{k}(r).
\end{eqnarray}

We will denote   $\Omega= {B(0,r+\epsilon) \setminus ( B(0,r-\epsilon) \cup B(z,\epsilon)})$ for the purpose of notational brevity .  Then, for an arbitrary but fixed $\epsilon$ one can interchange the order of integration (according to Fubini's Theorem) in the following way:
\begin{eqnarray}\label{eq2}
\nonumber c^\epsilon_{k}(r)&=&{i \over 4 \pi^2}  \int^{2 \pi}_0  \int  \int_{B(0,r-\epsilon)} h(\xi)  { e^{-i k \theta} \over (\xi -r e^{i \theta})^2} \  d \bar{\xi} \wedge d \xi d \theta\\
\nonumber &+&{i \over 4 \pi^2 }  \int^{2 \pi}_0  \int  \int_{B(0,R) \setminus B(0,r+\epsilon)} h(\xi) { e^{-i k \theta} \over (\xi -r e^{i \theta})^2} \  d \bar{\xi} \wedge d \xi d \theta\\
\nonumber &+&{i \over 4 \pi^2}  \int^{2 \pi}_0  \int  \int_{\Omega} h(\xi)  { e^{-i k \theta} \over (\xi -r e^{i \theta})^2} \   d \bar{\xi} \wedge d \xi d \theta\\
\nonumber &=&-{1 \over 2 \pi^2}  \int^{r-\epsilon}_0  \int^{2 \pi}_0  h(\rho e^{i \alpha}) \int^{2 \pi}_0  { e^{-i k \theta} \over (\rho e^{ i \alpha} -r e^{i \theta})^2} \  d \theta d \alpha  \rho d \rho \\ 
\nonumber &-&{1 \over 2 \pi^2}  \int^{R}_{r+\epsilon}  \int^{2 \pi}_0  h(\rho e^{i \alpha}) \int^{2 \pi}_0  { e^{-i k \theta} \over (\rho e^{ i \alpha} -r e^{i \theta})^2} \  d \theta d \alpha  \rho d \rho \\ 
&+&{i \over 4 \pi^2}  \int^{2 \pi}_0  \int  \int_{\Omega} h(\xi)  { e^{-i k \theta} \over (\xi -r e^{i \theta})^2} \  d \bar{\xi} \wedge d \xi d \theta.
\end{eqnarray}

Next, define  $\Omega_2= B(z,\epsilon^\beta) \setminus B(z,\epsilon)$ for some real positive $\beta<1$ to be specified later ($B(z,\epsilon) \subset B(z,\epsilon^\beta)$ since $\epsilon$ is always less than $1$ by our choice of $\delta$), and $\Omega_1=\Omega \setminus \Omega_2$. Using H\"older continuity once again, one can show that the difference  
\begin{eqnarray} \label{I_eps}
\quad \quad \quad I_\epsilon&=&\left|  {i \over 4 \pi^2} \int^{2 \pi}_0 \!\! \int  \int_{\Omega} { h(\xi) e^{-i k \theta} \over (\xi -z)^2} \ d \bar{\xi} \wedge d \xi d \theta-{i \over 4 \pi^2} \int^{2 \pi}_0  \!\!  \int  \int_{\Omega} { h(z) e^{-i k \theta} \over (\xi -z)^2} \  d \bar{\xi} \wedge d \xi d \theta \right|
\end{eqnarray}
can be bounded in the following way:
\begin{eqnarray}
\nonumber I_\epsilon & \le&  {1 \over 4 \pi^2} \int^{2 \pi}_0  \int  \int_{\Omega} { \left| h(\xi)-h(z) \right| \over |\xi -z|^2} \   d \bar{\xi} \wedge d \xi d \theta \le {1 \over 4 \pi^2} \int^{2 \pi}_0   \int  \int_{\Omega} A   |\xi -z|^{\gamma-2}  \  d \bar{\xi} \wedge d \xi d \theta \\
 \nonumber & \le&   {A \over 2 \pi}  \int  \int_{\Omega}  \rho^{\gamma-2} \  d \alpha \rho d \rho \le  {A \over 2 \pi} \int  \int_{\Omega_1}  \rho^{\gamma-2} \  d \alpha \rho d \rho   + {A \over 2 \pi} \int  \int_{\Omega_2}  \rho^{\gamma-2} \  d \alpha \rho d \rho  \\
 \nonumber & \le&    {A \over 2 \pi} \epsilon^{\beta(\gamma-2)}{\rm Area}(\Omega_1)+ {A \over \gamma} (\epsilon^{\beta \gamma}-\epsilon^\gamma) \le  {A \over 2 \pi} \epsilon^{\beta(\gamma-2)} 4 \pi r \epsilon+ {A \over \gamma} (\epsilon^{\beta \gamma}-\epsilon^\gamma)  \\
  \nonumber & \le& 2 A r  \epsilon^{\beta(\gamma-2)+1}+ {A \over \gamma} (\epsilon^{\beta \gamma}-\epsilon^\gamma).
\end{eqnarray}

Notice, that 
\begin{eqnarray}\label{I_eps_lim}
\lim_{\epsilon \rightarrow 0} I_\epsilon&=&0
\end{eqnarray}
 as long as $\beta (2-\gamma)<1$. Now, it is not hard to show that
\begin{equation}\label{Daripa_eqn}
{i \over 4 \pi^2 }  \int  \int_{\Omega} { d \bar{\xi} \wedge d \xi  \over (\xi -z)^2}={1 \over 2 \pi} {(r-\epsilon)^2 \over z^2}
\end{equation}
\noindent (see \cite{Daripa2}), and, therefore, combining (\ref{I_eps_lim}) and (\ref{Daripa_eqn}), that 
\begin{eqnarray}
\nonumber\left| \!\!\!\!\!\!\phantom{  \int^{2 \pi}_0  }\!\!\!\!\!\! \right.&{i \over 4 \pi^2}&\left. \lim_{\epsilon \to 0} \int^{2 \pi}_0  \!\!\! \int   \int_{\Omega} { h(\xi) e^{ -i k \theta} \over (\xi -z)^2} \ d \bar{\xi} \wedge d \xi d \theta -{1 \over 2 \pi}\int^{2 \pi}_0  h(z) e^{- i k \theta} {r^2 \over z^2} \ d \theta \right|=\\
\nonumber &=&\left|  {i \over 4 \pi^2} \lim_{\epsilon \to 0}  \int^{2 \pi}_0  \!\!\! \int  \int_{\Omega} { h(\xi) e^{-i k \theta} \over (\xi -z)^2} \ d \bar{\xi} \wedge d d \xi d \theta - {1 \over 2 \pi} \lim_{\epsilon \to 0}  \int^{2 \pi}_0  h(z) e^{-i k \theta}  {(r-\epsilon)^2 \over z^2} \ d \theta  \right|\\
\nonumber &=&\lim_{\epsilon \to 0} \left|  {i \over 4 \pi^2}  \int^{2 \pi}_0  \!\!\!  \int  \int_{\Omega} { h(\xi) e^{-i k \theta} \over (\xi -z)^2} \ d \bar{\xi} \wedge d \xi d \theta - {1 \over 2 \pi}  \int^{2 \pi}_0  h(z) e^{-i k \theta} {(r-\epsilon)^2 \over z^2} \ d \theta \right| \\
\nonumber &=&\lim_{\epsilon \to 0} \left|  {i \over 4 \pi^2}  \int^{2 \pi}_0  \!\!\!  \int  \int_{\Omega} { h(\xi) e^{-i k \theta} \over (\xi -z)^2} \ d \bar{\xi} \wedge d \xi d \theta -  {i \over 4 \pi^2}  \int^{2 \pi}_0   \!\!\! \int  \int_{\Omega} { h(z) e^{-i k \theta} \over (\xi -z)^2} \ d \bar{\xi} \wedge d \xi d \theta  \right| \\
\nonumber &  \le&  \lim_{\epsilon \to 0} I_\epsilon=0,
\end{eqnarray}
\noindent where in the first equality we again used the Lebesgue Dominated Convergence Theorem. We conclude that 
\begin{eqnarray}\label{reduce_int}
{i \over 4 \pi^2} \lim_{\epsilon \to 0} \int^{2 \pi}_0  \int   \int_{\Omega} { h(\xi) e^{ -i k \theta} \over (\xi -z)^2} \ d \bar{\xi} \wedge d \xi d \theta&=&{1 \over 2 \pi}\int^{2 \pi}_0  h(z) e^{- i k \theta} {r^2 \over z^2} \ d \theta.
\end{eqnarray}

Observe, that
  \begin{equation}\label{h_kplus2}
{1 \over 2 \pi} \int^{2 \pi}_0   h(z) e^{-i k \theta} {r^2 \over z^2} \ d \theta={1 \over 2 \pi} \int^{2 \pi}_0  { e^{-i k \theta} \sum^{\infty}_{n=-\infty}h_n(r) e^{i n \theta} \over e^{2 i \theta}} \ d \theta=h_{k+2}(r).
\end{equation} 

The claim of the Corollary follows from (\ref{eq2}), (\ref{reduce_int}) and (\ref{h_kplus2}).
 \end{proof}
\bigskip

Now, consider 
\begin{eqnarray}
\nonumber I(\xi,r)=\int^{2 \pi}_0 \!\! { e^{-i k \theta} \over (\xi -r e^{i \theta})^2} \ d \theta &=&{1 \over r^2} \int^{2 \pi}_0 \!\!  { e^{-i k \theta} \over (\xi / r  -e^{i \theta})^2} \  d \theta ={1 \over i r^2} \oint_{\field{T}^1} \!\!  { 1 \over w^{k+1}(\xi / r  -w)^2} \ d w. 
\end{eqnarray}

We will first compute this integral in the case when  $|\xi|/r>1$:
 \begin{eqnarray}
\nonumber  I(\xi,r)&=& {2 \pi \over r^2} {\rm Res} \left[ {1 \over w^{k+1}(\xi / r  -w)^2 } \right]_{w=0}=-B_k \pi {r^k \over \xi^{k+2}},
\end{eqnarray}
\noindent where 
\begin{eqnarray}
\nonumber  B_k&=&\left\{-2 (k+1),  \ k \ge 0, \atop 0, \ k<0. \right.
\end{eqnarray}

Similarly, for $|\xi|/r<1$,
 \begin{eqnarray}
\nonumber I(\xi,r) &=& {2 \pi \over r^2} {\rm Res} \left[ {1 \over w^{k+1}(\xi / r  -w)^2 } \right]_{w={\xi \over  r}}\!\!\!+{2 \pi \over r^2} {\rm Res} \left[{1 \over w^{k+1}(\xi / r  -w)^2 } \right]_{w=0}\!\!=-A_k \pi {r^{k} \over  \xi^{k+2}},
\end{eqnarray}
\noindent where 
\begin{eqnarray}
\nonumber A_k&=&\left\{ 0 \ , \ k \ge 0,  \atop  2  (k+1) \ , k<0. \right.
\end{eqnarray}

Finally, we obtain:
\begin{eqnarray}
\nonumber c_k(r)&=&{A_k \over 2 \pi}  \int^r_0 \! \int^{2 \pi}_0 \!\! h(\rho e^{i \alpha}){r^k \over \xi^{k+2}} \ \rho d \rho d \alpha +{B_k \over 2 \pi}  \int^R_r \! \int^{2 \pi}_0  \!\!h(\rho e^{i \alpha})  {r^k \over \xi^{k+2}} \ \rho d \rho d \alpha +h_{k+2}(r)\\
\nonumber &=&{A_k \over 2 \pi}  \int^r_0\!\!  \int^{2 \pi}_0  {  r^k  \sum^{\infty}_{-\infty}  h_n(\rho) e^{i n \alpha }\! \!\!\over \rho^{k+1} e^{i(k+2) \alpha} } \  d \rho d \alpha +{B_k \over  2 \pi} \int^R_r \!\! \int^{2 \pi}_0 { r^k \sum^{\infty}_{-\infty} h_n(\rho) e^{i n \alpha } \!\!\! \over \rho^{k+1} e^{i(k+2) \alpha} } \  d \rho d \alpha\\ 
\nonumber &+& h_{k+2}(r)=A_k \int^r_0 { r^k  \over \rho^{k+1} }  h_{k+2}(\rho)   d \rho+B_k \int^R_r  { r^k  \over \rho^{k+1} }  h_{k+2}(\rho)   d \rho +h_{k+2}(r).
\end{eqnarray}

\subsection{Case $z \in \partial B(0,R)$}\label{sect4}
We will briefly outline a counterpart of Corollary \ref{Lebesgue} for $|z|=R$.
\begin{cor}
  \begin{eqnarray}
\nonumber    c_{k}(R)&=&-{1 \over 2 \pi^2} \lim_{\epsilon \to 0}  \int^{R-\epsilon}_0  \int^{2 \pi}_0  h(\rho e^{i \alpha}) \int^{2 \pi}_0  { e^{-i k \theta} \over (\rho e^{ i \alpha} -R e^{i \theta})^2} \  d \theta d \alpha  \rho d \rho.
  \end{eqnarray} 
\end{cor}  
  \begin{proof}
     According to Lemma \ref{lemma1}, we can exchange the order of the integration and the limit:
    \begin{eqnarray}
      \nonumber c_k(R)&=&{i \over 4 \pi^2}\int^{2 \pi}_0 \lim_{\epsilon \to 0} \int  \int_{B(0,R) \setminus B(z,\epsilon)} h(\xi) { e^{-i k \theta} \over (\xi -R e^{i \theta})^2} \  d \bar{\xi} \wedge d \xi d \theta \\
      \nonumber &=&{i \over 4 \pi^2} \lim_{\epsilon \to 0} \int^{2 \pi}_0 \int  \int_{B(0,R) \setminus B(z,\epsilon)} h(\xi)  { e^{-i k \theta} \over (\xi -R e^{i \theta})^2} \  d \bar{\xi} \wedge d \xi d \theta= \lim_{\epsilon \to 0}  c^\epsilon_{k}(R).
    \end{eqnarray}

Denote  $\Omega= B(0,R) \setminus \left( B(0,R-\epsilon) \cup B(z,\epsilon) \right)$.  Then, we obtain for an arbitrary but fixed $\epsilon$:
\begin{eqnarray}\label{c_k_eps}
\nonumber c^\epsilon_{k}(R)&=&-{1 \over 2 \pi^2}  \int^{R-\epsilon}_0  \int^{2 \pi}_0  h(\rho e^{i \alpha}) \int^{2 \pi}_0  { e^{-i k \theta} \over (\rho e^{ i \alpha} -R e^{i \theta})^2} \  d \theta d \alpha  \rho d \rho \\ 
&+&{i \over 4 \pi^2}  \int^{2 \pi}_0  \int  \int_{\Omega} h(\xi)  { e^{-i k \theta} \over (\xi -R e^{i \theta})^2} \  d \bar{\xi} \wedge d \xi d \theta.
\end{eqnarray}

Again, define  $\Omega_2= B(z,\epsilon^\beta) \setminus B(z,\epsilon)$, for some real positive $\beta<1$, and $\Omega_1=\Omega \setminus \Omega_2$. As in Corollary \ref{Lebesgue}, one shows that the difference (\ref{I_eps}) can be bounded  by $A r  \epsilon^{\beta(\gamma-2)+1}+ {(A / 2 \gamma}) (\epsilon^{\beta \gamma}-\epsilon^\gamma)$, and that $\lim_{\epsilon \rightarrow 0} I_\epsilon=0$ as long as $\beta (2-\gamma)<1$. However, $h(z) \equiv 0$ for all $z \in \partial B(0,R)$, hence
\begin{eqnarray}
\nonumber I_\epsilon & =&\left|  {i \over 4 \pi^2} \int^{2 \pi}_0  \int  \int_{\Omega} { h(\xi) e^{-i k \theta} \over (\xi -z)^2} \ d \bar{\xi} \wedge d \xi d \theta \right|
\end{eqnarray}
and therefore
\begin{equation} \label{eqn_zero}
\nonumber {i \over 4 \pi^2} \lim_{\epsilon \to 0} \int^{2 \pi}_0  \int  \int_{\Omega} h(\xi)  { e^{-i k \theta} \over (\xi -R e^{i \theta})^2} \  d \bar{\xi} \wedge d \xi d \theta=0.
\end{equation}

The claim of the Lemma follows from  (\ref{c_k_eps}) and (\ref{eqn_zero}). 
\end{proof} 

\bigskip

To obtain an expression for $c_k(R)$, consider the integral $\int^{2 \pi}_0  { e^{-i k \theta} / (\xi -R e^{i \theta})^2} \ d \theta$. It is not hard to show that this integral is equal to $-A_k \pi  {R^k / \xi^{k+2}}$,  and, therefore,
\begin{eqnarray}
\nonumber c_k(R)&=&-{1 \over 2 \pi^2}  \lim_{\epsilon \to 0} \int^{R-\epsilon}_0  \int^{2 \pi}_0  h(\rho e^{i \alpha}) \int^{2 \pi}_0  { e^{-i k \theta} \over (\rho e^{ i \alpha} -R e^{i \theta})^2} \  d \theta d \alpha  \rho d \rho \\
\nonumber &=&{1 \over 2 \pi}  \int^{R}_0  \int^{2 \pi}_0  h(\rho e^{i \alpha}) A_k  {R^k \over \xi^{k+2}} \ d \alpha  \rho d \rho \\
\nonumber &=&{A_k \over 2 \pi}  \int^{R}_0  \int^{2 \pi}_0   { R^k  \sum^{\infty}_{n=-\infty} h_n(\rho) e^{i n \theta } \over \rho^{k+1} e^{i(k+2) \theta} } \  d \rho d \theta \\
\nonumber  &=&A_k \int^{R}_0  { R^k  \over \rho^{k+1} }  h_{k+2}(\rho)  \  d \rho.
\end{eqnarray}

\subsection{Case $z=0$} \label{sect3}
We will prove a counterpart of Corollary (\ref{Lebesgue}) at the same time obtaining an expression for $c_k(0)$.
\begin{cor} \label{Lebesgue2}
\begin{eqnarray}
\nonumber c_k(0)&=&\left\{ 
\begin{array}{cc} 
	0,&  k \ne 0,\\
	-2  \int^{R}_0 {h_2(\rho) \over \rho }  d \rho ,  & k=0.
      \end{array}
\right.
\end{eqnarray}
\end{cor}
\begin{proof} 
As before, according to Lemma \ref{lemma1} for an arbitrary but fixed $\epsilon$ one can use the Lebesgue Dominated Convergence Theorem to change the order of the integration and the limit:
\begin{eqnarray}
\nonumber c_{k}(0)&=&{i \over 4 \pi^2}  \int^{2 \pi}_0  \lim_{\epsilon \to 0}  \int  \int_{B(0,R)\setminus B(0,\epsilon)} h(\xi)  { e^{-i k \theta} \over \xi^2} \  d \bar{\xi} \wedge d \xi d \theta\\
\nonumber &=&{i \over 4 \pi^2}  \lim_{\epsilon \to 0} \int^{2 \pi}_0   \int  \int_{B(0,R)\setminus B(0,\epsilon)} h(\xi)  { e^{-i k \theta} \over \xi^2} \  d \bar{\xi} \wedge d \xi d \theta\\
\nonumber &=&-{1 \over 2 \pi^2}  \lim_{\epsilon \to 0}  \int^{R}_\epsilon  \int^{2 \pi}_0  {h(\rho e^{i \alpha})   \over \rho^2 e^{ 2 i \alpha}} \int^{2 \pi}_0  e^{-i k \theta} \  d \theta d \alpha  \rho d \rho \\ 
\nonumber &=&-{1 \over 2 \pi^2}  \int^{R}_0 \int^{2 \pi}_0  {\sum^{n=\infty}_{n=-\infty} h_n(\rho) e^{i n \alpha}  \over \rho^2 e^{ 2 i \alpha}} \int^{2 \pi}_0  e^{-i k \theta} \  d \theta d \alpha  \rho d \rho \\ 
\nonumber &=&\left\{ 
\begin{array}{cc} 
	0,&  k \ne 0,\\
	-2  \int^{R}_0 {h_2(\rho) \over \rho }  d \rho ,  & k=0.
      \end{array}
\right.
\end{eqnarray}

The H\"older property of $h$ and the fact that $h_2(0)=0$ show that the last integral ($k=0$) is convergent.
\end{proof}

\subsection{Case $z \in \field{C} \setminus \overline{B (0,R)}$}\label{sect5}
In this case the limiting procedure and the angular integration in $(\ref{ck_long})$  can be interchanged right away. 
\begin{eqnarray}
\nonumber \int^{2 \pi}_0  { e^{-i k \theta} \over (\xi -r e^{i \theta})^2} \ d \theta= {2 \pi \over r^2} {\rm Res} \left[ {1 \over  w^{k+1}(\xi / r -w)^2 } \right]_{w=0}=-A_k \pi {r^k \over \xi^{k+2}},
\end{eqnarray}
\noindent and
\begin{eqnarray}
\nonumber c_k(r)&=&{i \over 4 \pi^2}  \int  \int_{B(0,R)} h(\xi) A_k {r^k \over \xi^{k+2}} \  d \bar{\xi} \wedge d \xi \\
\nonumber & =& {A_k \over 2 \pi}  \int  \int_{B(0,R)} \sum^{\infty}_{n=-\infty} h_n(\rho) e^{i n \theta } {r^k \over \xi^{k+2}} \ \rho d \rho d \theta \\
\nonumber  &=&A_k \int^R_0  { r^k  \over \rho^{k+1} }  h_{k+2}(\rho)  \  d \rho.
\end{eqnarray}

\bigskip

\section{Recursive formulae}\label{sect6}
\setcounter{equation}{0}
Define
\begin{eqnarray}
\nonumber d_k(r,r')&=& 2 (k+1) \int^{r'}_{r} {1 \over \rho} \left({ r'  \over \rho } \right)^k  h_{k+2}(\rho)  \  d \rho, 
\end{eqnarray}
then, using (\ref{Hilbert_alg_1})--(\ref{A_k_B_k}), one readily gets the following recursive formula for the coefficients of the Hilbert transform:
\begin{eqnarray}\label{H_rec_formula}
    c_k(r)=\left(r \over r'  \right)^k (c_k(r')-h_{k+2}(r')-d_k(r,r')) + h_{k+2}(r), 
\end{eqnarray}
which holds for any pair of nonnegative radii $r$ and $r'$ for which the ratio $(r/r')^k$ is defined.

We will also note, that it is convenient to compute $p_k(r)$ in (\ref{Cauchy_alg}) using the following recursive
relations:
\begin{eqnarray}\label{P_rec_formula}
 p_k(r)&=&\left(\frac{r}{r'}\right)^k p_k(r')-e_k(r,r'),\ e_k(r,r')=2\int_r^{r'}\left(\frac{r}{\rho}\right)^k h_{k+1}(\rho)d\rho.
\end{eqnarray}

 These recursive formulae are extremely useful because they allow to compute the value of a coefficient at the $i$-th radius using its value at either $(i+1)$-th, or $(i-1)$-th radius, thus avoiding the full integration in formulae  (\ref{Cauchy_alg}) and (\ref{Hilbert_alg_1})--(\ref{A_k_B_k}).

For reference purposes we will rewrite formulae  (\ref{H_rec_formula}) and (\ref{P_rec_formula}) for a piecewise-constant  extrapolation of data $h_k(r_i)$ in Appendix A. Strictly speaking, our derivation of the expressions for the Fourier coefficients of the Hilbert transform does not hold for a discontinuous $h$ with piecewise-constant coefficients. However, given an $\epsilon>0$ there always is a H\"older-continous function whose coefficients are $\epsilon$-close to $h_k(r)$ in the $L_1$ norm, and, therefore, there always is a H\"older-continous smoothening of piecewise-constant data such that the coefficients of its Hilbert transform are arbitrarily close to those given by expressions (\ref{Hilbert_1})--(\ref{Hilbert_3}) in the sup norm.

\section{Another way to compute $T$}
\setcounter{equation}{0}
In this Section we propose an alternative way for computing $T [h]$. By a theorem, presented in ~\cite{Ahlfors}, the following identity holds for all $h\in{}C^2_0$ (twice differentiable $h$ with a compact support):
\begin{equation}  \label{eq:(Ph)_w}
  T[h](z)=\partial_z P[h] (z)=\partial_z \P[h](z).
\end{equation}

As we have already mentioned,  in \cite{Daripa1} Daripa has suggested an effective algorithm for a numerical computation of integral transform (\ref{eq:Pt}). The algorithm in \cite{Daripa1} for computing $\P$ is defined in
terms of a circular $N \times M$ grid and its computational complexity
is $O(N M\log_2{M})$. In order to use this algorithm it is necessary to
find an approximation for $h_z$, given the values of $h$ on a circular
grid. A formula for $h_z$ in polar coordinates is supplied in the following Lemma, whose easy prove will be omitted.

\begin{lemma}
  Suppose that $f$ is $C^1(\R^2)$, and in polar coordinates $f$ is represented by $g$, i.e. $f(x,y)=g(r,\th)$, where $x+iy=re^{i\th}$. Then
  \begin{eqnarray}
    \label{eq:f_z_in_g_r} f_z(x,y)&=&\frac{e^{-i\th}}{2}(g_r(r,\th)+g_\th(r,\th)/ir),\\
    \label{eq:f_zbar_in_g_r} f_\z(x,y)&=&\frac{ie^{i\th}}{2}(-ig_r(r,\th)+g_\th(r,\th)/r).
  \end{eqnarray}
\end{lemma}

\section{Comparison of the two algorithms}

We have implemented the algorithm for the Cauchy transform (see formulae (\ref{eq:Pt})--(\ref{Cauchy_alg}) and (\ref{P_rec_formula})  ), both algorithms for the Hilbert transform ((\ref{T_coeff})--(\ref{A_k_B_k}), (\ref{H_rec_formula}) and (\ref{eq:(Ph)_w})) for the piecewise linear approximation of Fourier coefficients as a set of routines in the programming language Ada 95 (cf \cite{ADA} for the language standard). We have used the public version 3.15p of the GNAT compiler \cite{GNAT}. Our programs can be found at \cite{programs}.

 Of particular interest to us will be a comparison of the two implementations of the iteration scheme  $h \mapsto  T[\mu(h+1)]$: in one of them, the Hilbert transform has been computed through formula (\ref{H_rec_formula}), in the other, as the $z$-derivative of transform (\ref{eq:Pt}). Below, we will refer to these two iteration schemes  as Scheme 1 and Scheme 2. 

The prototypical Beltrami differential for algorithm testing in the literature seems to be $\mu(z)=z^p \bar{z}^q$ in some disk and $\mu(z)=0$ outside of that disk. Both Hilbert and Cauchy transforms of this differential, as well as the solution of the Beltrami equation itself,  can be found exactly. Here we will quote the formula for the Hilbert transform:
\begin{eqnarray}
\nonumber T[z^p \bar{z}^{q}]&=& { p \ \! \eta(R-|z|)  \over q+1} z^{p-1} \bar{z}^{q+1}+\\
\label{T_zpzq} & + & \left[\eta(q+1-p)- \eta(R-|z|) \right] {p-q-1 \over q+1} R^{2(q+1)} z^{p-q-2}\!\!,
\end{eqnarray}
where $\eta$ is the Heaviside step function. The derivation of this formula is relatively straightforward, and can be done either by a direct calculation or by differentiating the expression for the Cauchy transform of $z^p \bar{z}^{q}$ which have appeared in several publications (cf e.g. \cite{Daripa3}).

This formula suggests yet one more algorithm for the Hilbert transform of $C^k_0$ functions representable as finite power series on its compact support. Let $h(z,\bar{z})=\eta(R-|z|) \sum^{m,n}_{i,j=0} a_{i,j} z^i \bar{z}^j$ be such function. Then one can first use the linearity of the Hilbert transform to represent $T[h]$ as $T[h](z,\bar{z})=\sum^{m,n}_{i,j=0} a_{i,j} T[z^i \bar{z}^j]$, compute each $T[z^i \bar{z}^j]$ in this series using (\ref{T_zpzq}) and sum up. We will refer to this procedure as Scheme 3. Of course, this Scheme has a limited value, as many differentials do not have sufficient smoothness, or are sometimes given in numerical experiments simply as an array of values on a grid.


Here, we will use a quartic test Beltrami differential which readily yields itself to the procedure that we just described:
\begin{eqnarray} \label{test_diff}
  \mu(z)&=& \left\{ 
\begin{array}{cc}
   a-{a \over s^2} z \bar{z} +{a  \over 4 s^4} z^2 \bar{z}^2, & |z| \le \sqrt{2}  s,\\
   0, &   |z| > \sqrt{2} s,
\end{array}
\right.
\end{eqnarray}
where $a$ and $s$ are some real positive parameters. Our $\mu$ is given by a real-valued, non-negative, 
radially symmetric $C^2_0$ function on $\field{C}$ with a single maximum at zero, compactly supported in the disk 
of radius $\sqrt{2} s$. We would like to emphasize, that iterations $h^{n+1} \mapsto T[h^n (\mu+1)]$ can 
be performed {\it exactly} using Scheme 3 (that is, the only error is the rounding error of our computing 
device). In particular, there is no need for any approximation of the functions being transformed, as it 
was the case with the choice of the spline approximation of the coefficients $h_k$ in formulae (\ref{Hilbert_alg_1})--(\ref{Hilbert_alg_2}). We will use Scheme 3 as a measuring stick for the accuracy of Schemes 1 and 2. 

Let $g(r,\theta)=\sum_k g_k(r) e^{i k \theta}$ be smooth, then it follows from (\ref{eq:f_z_in_g_r}) that the Fourier coefficients of the $z$-derivative of this function, denoted by $\tilde{g}$, are given by the following exact expression:
\begin{eqnarray}
  \tilde{g}_{k-1}(r)&=&{g'_k(r)+k g_k(r)/r   \over 2}.
\end{eqnarray}

We have used this formula to compute the coefficients of the $z$-derivative of the iteration $h \mapsto \partial_z  P[ \mu(h+1)]$ for our choice (\ref{test_diff}) of $\mu$; $g'_k(r)$'s were supplied by right two-point approximations through values of $g_k(r)$ on the radial grid. 

We should also mentioned that we have used the fast Fourier transform (FFT) together with the point-wise multiplication of function values on the grid  to perform the multiplication of $\mu$ and $h$ in both schemes. Indeed, the computational cost of the FFT is of order $O( N M \log_2 M )$, that of the point-wise multiplication of functional values is $O( N  M)$; at the same time, the complexity of the multiplication of two Fourier series representing $\mu$ and $h$ would be $O( N M^2)$.

We subjected the two schemes to several test on an 1.5 GHz Intel Pentium class PC. Our first test was that of convergence of the two iteration schemes. Table 1 compares the convergence rate (measured by the norm $\|h^n-h^{n-1} \|_\infty$ rounded to 3 significant figures) for both iteration schemes and several values of $n$. We have used zero initial conditions, $h^0=\mu$ in both cases. We have fixed the values of the parameters in the Beltrami differential to be $a=0.5$ and $s=1.0$. The dimensions of the radial grid were $N=500 $ (i.e. the disk of radius $\sqrt{2}$ contained $500$ radii of the grid) and $M=512$.  

\bigskip

\begin{center}
\begin{tabular}{|c|c|c|c|c|}
\hline
number of iterations, n $\phantom{\int^N_N \!\!\!\!\!\!\!\!\!\!\!}$ & 5 & 10 & 20 & 50 \\
\hline
Scheme 1 $\phantom{\int^N_N \!\!\!\!\!\!\!\!\!\!\!}$ & $1.17 \times 10^{-6}\!$ & $3.83 \times 10^{-11}\!$ & $1.80 \times 10^{-14}\!$  & $1.17 \times 10^{-19}\! $  \\ 
\hline
Scheme 2  $\phantom{\int^N_N \!\!\!\!\!\!\!\!\!\!\!}$ & $1.20 \times 10^{-6}\!$ & $3.91 \times 10^{-10}\!$ & $6.79 \times 10^{-14}\!$  & $8.27 \times 10^{-20}\! $  \\ 
\hline
\end{tabular}
\end{center}
\bigskip
\begin{center}
Table 1.   Convergence rate ($\|h^n-h^{n-1} \|_\infty$) for the two iteration schemes. 
\end{center}

\bigskip

It is not surprising that the convergence rates of the schemes are comparable: These rates should depend only on the value of the essential supremum of $\mu$ and not on the particular realization of the Hilbert transform. 

Next, we have measured the execution time (rounded to nearest seconds) of $10$ iteration steps for several grids. Table 2 displays the results (the values of $a$, $s$ and $h^0$ were as before).

\bigskip

\begin{center}
\begin{tabular}{|c|c|c|c|c|}
\hline
$N \times M$ $\phantom{\int^N_N \!\!\!\!\!}$ &  $500 \times 256$ & $1000 \times 512$ & $5000 \times 1024$ & $7000 \times 2048$ \\
\hline
Scheme 1  $\phantom{\int^N_N \!\!\!\!\!}$ & 8  & 31 & 313 & 965 \\ 
\hline
Scheme 2  $\phantom{\int^N_N \!\!\!\!\!}$ & 10 & 38 & 375 & 1166 \\ 
\hline
\end{tabular}
\end{center}

\bigskip

\begin{center}
Table 2.   Execution time for several grids (in sec). 
\end{center}

\bigskip

As it should be expected, Scheme 2 takes somewhat longer: The computational complexity of both integral transforms is the same, $O(N M \log_2{M})$, however, the iteration $h \mapsto  \partial_z P[ \mu(h+1)]$ also involves a differentiation whose computational complexity is $O(N M)$. Nevertheless, for the tested range of grids, the speed benefits of Scheme 1 probably do not justify an extra effort in programming the Hilbert transform: As we have already mentioned, in our experience writing and testing such code might be rather cumbersome and lengthy, and could significantly increase the project development time.

The $l_\infty$ norm of the difference of the outcomes of Schemes 1 and 3 and Schemes 2 and 3 after 
10 steps (with the 
same parameters and initial conditions as before)  for several grids is reported in Table 3. As we have 
already explained it, Scheme 3 serves as a measuring stick for our purposes. 

\bigskip

\begin{center}
\begin{tabular}{|c|c|c|c|c|}
\hline
$\phantom{\int^N_N \!\!\!\!\!} N \times M$             &  $500 \times 256$       & $1000 \times 512$      & $5000 \times 1024$    & $7000 \times 2048$ \\
\hline 
$\phantom{\int^N_N \!\!\!\!\!}\|h^{2}_1-h^{2}_3\|_{\infty}$  & $4.00 \times 10^{-7}$   & $1.00 \times 10^{-7}$ & $7.12 \times 10^{-8}$    & $1.89  \times 10^{-7}$  \\ 
\hline
$\phantom{\int^N_N \!\!\!\!\!}\|h^{2}_2-h^{2}_3 \|_{\infty}$ & $8.39 \times 10^{-5}$  &  $4.20 \times 
10^{-5}$ & $8.39 \times 10^{-6}$ & $5.99 \times 10^{-6}$ \\ 
\hline
\end{tabular}
\end{center}
\bigskip
\begin{center}
Table 3.   Comparison of the performance of Schemes 1 and 2 against Scheme 3. 
\end{center}

\bigskip

With regards to the accuracy of Schemes 1 and 2,  Table 3 demonstrates that Scheme 1 is the scheme of choice. A bigger error in Scheme 2 can be attributed to the fact that we are using a numerical approximation for the differentiation  in the iteration $h \mapsto  \partial_z P[ \mu(h+1)]$. This approximation introduces an extra error which is absent in Scheme 1. However, if the degree of differentiability of $h$ is known, one can use a proper spline approximation of the derivative and thus significantly reduce this error. We still feel that in those cases when a high accuracy is not required, Scheme 2 can be more preferable.

\section{Appendix A}
We will now return to formula (\ref{H_rec_formula}). Consider a collection of real numbers  $0=r_1 < r_2 < \ldots < r_{L-1} <  r_L=R < r_{L+1} < \ldots < r_N$ that specify a radial grid.  We will denote the midpoint of the interval $[r_i,r_{i+1}]$ by $\varrho_i$; $\varrho_{0}$ will be $0$ by definition.
 
If  $h_k$'s are compactly supported {\it piecewise-constant} functions
\begin{eqnarray}
  \nonumber h_k(r)&=& \left\{
  \begin{array}{cc} 
    h_{k,i},  & r  \in [\varrho_{i-1},\varrho_i), 1 \le i \le L,  \\
      0,    & r > \varrho_{L+1},
  \end{array}
  \right.
\end{eqnarray}
then the coefficients of the Cauchy and the Hilbert transforms are given by the following set of expressions.

\bigskip

\noindent{\it Coefficients of the Hilbert transform for  piecewise-constant data}.

\noindent $\bullet$ $r=r_1$
\begin{eqnarray}\label{Hilbert_1}
c_k(r_1)&=&\left\{ 
\begin{array}{cc} 
  c_0(r_2)-s_2(r_2) -2 s_2(r_2) \ln{r_2 \over \varrho_0 }, &  k=0, \\
   0,&  k \ne 0.
\end{array}
  \right.
\end{eqnarray}
\noindent $\bullet$ $r \in (\varrho_{i-1},\varrho_i]$, $1 \le i \le L$,
  \begin{eqnarray}\label{Hilbert_2}
    c_k(r)&=&\!\left\{ \!\!\!
    \begin{array}{cc} 
      \left( {r \over r_{i+1}} \right)^k \left(c_k(r_{i+1})-s_{k+2}(r_{i+1})\right)-d_{k,i}(r) +s_{k+2}(r), & \! k \ge 1,\\
      \! c_0(r_{i+1})\!-\!s_2(r_{i+1}) \!-\!2 s_2(r_{i+1}) \ln{r_{i+1} \over \varrho_i }\!-\!2 s_2(r_i) \ln{\varrho_i \over r }\!+\!s_2(r), & \! k=0, \\
      \left( {r \over r_{i-1}} \right)^k \left(c_k(r_{i-1})-s_{k+2}(r_{i-1})\right)+b_{k,i}(r) +s_{k+2}(r),& \! k < 0,
    \end{array}\!\!\!\!\!\!\!\!\!\!\!\!
    \right. 
  \end{eqnarray}
  where 
  \begin{eqnarray}
  \nonumber d_{k,i}(r)&=&2 {k+1 \over k} \left[s_{k+2}(r_{i+1}) \left({r^k \over \varrho^k_i}-{r^k \over r^k_{i+1}} \ \right) + s_{k+2}(r_i) \left(1-{r^k \over \varrho^k_i } \right) \right], \\
  \nonumber b_{k,i}(r)&=&2 {k+1 \over k} \left[s_{k+2}(r_{i-1}) \left({r^k \over r^k_{i-1}}-{r^k \over \varrho^k_{i-1}} \right) + s_{k+2}(r_i) \left({r^k \over \varrho^k_{i-1} }-1\right) \right].
\end{eqnarray}
\noindent$\bullet$ $r \in (\varrho_{i-1},\varrho_i]$, $L < i \le N$,
  \begin{eqnarray}\label{Hilbert_3}
      c_k(r)&=&\left\{ 
      \begin{array}{cc} 
	0,&  k \ne 0,\\
	{r^k \over r^k_{i-1}} c_k (r_{i-1}), &  k < 0.
      \end{array}
      \right. 
    \end{eqnarray}
    
    \bigskip

\noindent{\it Coefficients of the Cauchy transform for  piecewise-constant data}.

\noindent $\bullet$ $r=r_1$
\begin{eqnarray}	
  \nonumber c_k(r_1)&=&\left\{		 
  \begin{array}{cc} 		
    c_0(r_2)-2 s_1(r_2) \varrho_0 -2 s_1(r_2) (r_2-\varrho_1), &  k=0, \\
    0,&  k \ne 0.
  \end{array}	
  \right.
\end{eqnarray}
\noindent $\bullet$ $r \in (\varrho_{i-1},\varrho_i]$, $1 \le i \le L$,
  \begin{eqnarray}
    \nonumber c_k(r)&=&\left\{ 
    \begin{array}{cc} 
      \left( {r \over r_{i+1}} \right)^k c_k(r_{i+1})-d_{k,i}(r), &  k > 1,\\
	   {r \over r_{i+1}} c_1(r_{i+1}) -2 r s_2(r_i) \ln{\varrho_i \over r}-2 r s_2(r_{i+1})\ln{r_{i+1} \over \varrho_i} , &  k=1, \\
	   \left( {r \over r_{i-1}} \right)^k c_k(r_{i-1})+b_{k,i}(r),&  k <= 0,
    \end{array}
    \right. 
  \end{eqnarray}
  where 
  \begin{eqnarray}
    \nonumber d_{k,i}(r)&=&{2 r \over 1-k} \left[s_{k+1}(r_i) \left({r^{k-1} \over \varrho^{k-1}_i}-1 \right) + s_{k+1}(r_{i+1}) \left({r^{k-1} \over r^{k-1}_{i+1} }-{r^{k-1} \over \varrho^{k-1}_i } \right) \right], \\
    \nonumber b_{k,i}(r)&=&{2 r \over 1-k} \left[s_{k+1}(r_{i-1}) \left({\varrho^{1-k}_{i-1} \over r^{1-k}}-{r^{1-k}_{i-1} \over r^{1-k}} \right) + s_{k+1}(r_i) \left(1-{\varrho^{1-k}_{i-1} \over r^{1-k}}\right) \right].
  \end{eqnarray}
  \noindent $\bullet$ $r \in (\varrho_{i-1},\varrho_i]$, $L < i \le N$,
    \begin{eqnarray}
      \nonumber c_k(r)&=&\left\{ 
      \begin{array}{cc} 
	0,&  k \ne 0,\\
	{r^k \over r^k_{i-1}} c_k (r_{i-1}), &  k < 0.
      \end{array}
      \right. 
    \end{eqnarray}
  
\section{Acknowledgments} 

The author would like to thank Michael Yampolsky for his numerous helpful comments.

\end{document}